\documentclass[a4paper,12pt]{article}

\usepackage{color}
\usepackage[pdftex]{graphicx}
\usepackage{amsmath,amssymb,amsthm}
\usepackage{wrapfig}
\usepackage{ifpdf}

\newtheorem{theorem}{Theorem}

\newcommand{\dist}{\mathop{\rm dist}\nolimits}

\newcommand{\rp}{\mathop{\rm dim^E_2}\nolimits}
\newcommand{\sr}{\mathop{\rm  dim^S_2}\nolimits}

\newcommand{\cl}{\mathop{\rm Clq}\nolimits}

\title {Borsuk's conjecture for two--distance sets and its equivalent formulation for graphs}

\author {Oleg R. Musin}

\begin{document}
\date{}
\maketitle

\begin{abstract} 
Every graph $G$ can be embedded in a Euclidean space as a two--distance set. This allows us to reformulate the analogue of Borsuk's conjecture for two--distance sets in terms of graphs. This conjecture remains open for dimensions from 4 to 63. This short note also discusses an approach for finding counterexamples using graphs, as well as its generalization for $s$--distance sets. 
\end {abstract}

\section{Introduction}

Let $S$ be a bounded subset of $\mathbb R^d$ and $B(S)$ denote the smallest number of subsets with diameters smaller than the diameter of $S$ that cover $S$. Karol Borsuk conjectured that $B(S)\le d+1$. It is true for $d=2, 3$.  In 1993 Jeff Kahn and Gil Kalai showed that if $d\ge 2015$  then the conjecture is not true. In subsequent years, the dimension $d$ has been consistently reduced. In 2013 Andriy Bondarenko found a counterexample to Borsuk's conjecture for $d=65$ \cite{Bond14}. Soon after, Thomas Jenrich \cite{Jen}  discovered a 64-dimensional subset in Bondarenko's counterexample that also fails to satisfy Borsuk's conjecture. 

Bondarenko's counterexamples were obtained for two--distance sets.  A set $S$ in Euclidean space ${\Bbb R}^d$ is called a {\it two--distance set}, if there are two distances $a$ and $b$, and the
distances between pairs of distinct points of $S$ are either $a$ or $b$. We checked the case $d = 4$ and found no counterexample to Borsuk's conjectures for two--distance sets. This was later confirmed in \cite{Sz2018} by Ferenc Sz\"oll\"osi who enumerated all two--distance sets in four dimensions. Danylo Radchenko \cite{Rad} informed the author that he had not found such counterexamples for dimensions $d\le7$, his algorithm relies on  \cite{lisonek1997}.  

\medskip

Let $G=(V,E)$ be a simple graph,  $n=|V|$,  
$$
x_{ij}:=\left\{
\begin{array}{l}
1 \; \mbox{ if } e_{ij} \in E \\  %\mbox{ i.e. is an edge of } G\\
t \;  \mbox{ otherwise}
\end{array} \right. , \quad 
A_G(t):=
 \begin{bmatrix}
 0 & 1 & 1 & ... & 1\\
 1 & 0 & x_{12} & ... & x_{1n}\\
  1  & x_{21} & 0 & ... & x_{2n}\\
\hdotsfor{5}\\
\hdotsfor{5}\\
  1 & x_{n1} & x_{n2} & ... & 0
  \end{bmatrix}, \quad C_G(t)=\det(A_G(t)).
\eqno (1)
 $$
 
 %Suppose there is a solution of the algebraic equation $C_G(t)=0$ with $t>1$. Then denote by $\mu(G)$ the multiplicity of the smallest root $t>1$ of $C_G(t)=0$. If for all roots  $C_G(t)=0$ we have $t\le 1$, then we assume that $\mu(G)=0$.
 Suppose that there is a solution to the algebraic equation $C_G(t)=0$ for $t>1$. Then we denote by $\mu(G)$ the multiplicity of the smallest zero $t>1$ of this equation. If $C_G(t)$ has no zeros $t>1$, then we assume that $\mu(G)=0$.
 
 Let $\theta(G)$ denote the {\em clique covering number}, i.e. the smallest number of cliques that together include every vertex of  $G$.
 
 \begin{theorem} \label{main} Let $G$ be a simple graph with $n$ vertices.  Then $G$ can be embedded in $\mathbb R^d$ as a two--distance set $S$ with $d=n-\mu(G)-1$ and $B(S)=\theta(G)$. In particular, if  
   $$\theta(G)+\mu(G)> n$$ then $S$ is a counterexample for Borsuk's conjecture.  
\end{theorem}

In Section 2, we go into more detail about representing graphs as two--distance sets and prove this theorem. Much more details and examples are contained in our paper \cite{musin19}. 
 Section 3 discusses possible applications of Theorem 1 and its generalization for $s$--distance sets to finding counterexamples to Borsuk's conjecture using discrete optimization and machine learning methods.

%where $n=|V|$ and 
 %$$ x_{ij}:=\left\{
%\begin{array}{l}
%1 \; \mbox{ if } e_{ij} \in E,  \mbox{ i.e. is an edge of } G\\
%t \;  \mbox{ otherwise}
%\end{array} \right.
%$$

 \section{Graphs and two--distance sets} 
 \subsection{Two--distance sets.} 
Let $S=\{p_1,...,p_n\}$ be a two--distance in ${\Bbb R}^d$ with $a=1$ and $b>1$. Let  $d_{ij}:=\dist(p_i,p_j)$.
Consider the Cayley--Menger determinant 
 $$
C_S:=
 \begin{vmatrix}
 0 & 1 & 1 & ... & 1\\
 1 & 0 & d_{12}^2 & ... & d_{1n}^2\\
  1  & d_{21}^2 & 0 & ... & d_{2n}^2\\
\hdotsfor{5}\\
\hdotsfor{5}\\
  1 & d_{n1}^2 & d_{n2}^2 & ... & 0
  \end{vmatrix} 
%\eqno (2.2)
 $$
  Since for $i\ne j$, $d_{ij}=1$ or $b$,  $C_S$ is a polynomial in $t=b^2$. Here the connection between the two polynomials $C_G$ and $C_S$ is obvious. If $G$ is a graph with vertices $\{p_i\}$ and edges $e_{ij}$, where $d_{ij}=1$, then $C_G(t)=C_S(t)$. 
  
  %Denote this polynomial by $C_G(t)$. 

% Actually, in \cite{ES} instead of $C_G$ the discriminating polynomial $D(t)$ is considered. This polynomial can be defined through the Gram determinant. Since, see \cite[Lemma 9.7.3.3]{Berger},  
 %$$C_G(t)=(-1)^nD(t)$$  these two polynomials are the same up to the sign and therefore have the same roots. 

\medskip

Einhorn and Schoenberg \cite{ES} showed that  $d$ can be  explicitly expressed in terms of the multiplicity  $\mu(S)$ of the smallest root $t>1$ of  the equation $C_S(t)=0$:
$$d=n-\mu(S)-1.$$
Equivalently, see \cite[Sect. 2]{musin19}, for graphs we have 
$$ \rp(G)=n-\mu(G)-1, \eqno (2)$$
where $\rp(G)$ denotes the smallest dimension of the embedding of $G$ as a two--distance set.

  \subsection{Proof of Theorem \ref{main}} 
  
 \begin{proof} By (2) we can embed $G$ into $\mathbb R^d$, $d=n-\mu(G)-1$, as a two-distance set $S$. Since $B(S)=\theta(G)$, then $B(S)>d+1$ if and only if \, $\theta(G)+\mu(G)> n$.
  \end{proof}

 \subsection{Spherical two--distance sets.} 
In \cite{musin19} we consider representations of a graph $G$ as spherical two--distance sets.  Let $f$ be a Euclidean representation of  $G$ in ${\Bbb R}^d$ with minimum distance $a=1$.   We say that $f$ is {\em spherical} if the image $f(G)$ lies on a $(d-1)$--sphere in ${\Bbb R}^d$.  We denote by $\sr(G)$ the smallest $d$ such that $G$ is spherically representable in  ${\Bbb R}^d$.  

If $d\le n-2$, then  $f$ is uniquely defined up to isometry \cite[Sect. 2]{musin19}. Therefore, if $f$ is spherical, then the circumradius of $f(G)$ is also uniquely defined.  We denote it by ${\mathcal R}(G)$. If $f$ is not spherical or $\mu(G)=0$, then we put ${\mathcal R}(G)=\infty$.  By Lemma 3.1 in \cite{musin19} we have 
$${\mathcal R}(G)=\sqrt{F_G(\tau_1)},$$ 
where 
$$
F_G(t):=-\frac{1}{2}\frac{M_G(t)}{C_G(t)}, \qquad M_G(t):=
 \begin{vmatrix}
  0 & x_{12} & ... & x_{1n}\\
 x_{21} & 0 & ... & x_{2n}\\
\hdotsfor{4}\\
\hdotsfor{4}\\
 x_{n1} & x_{n2} & ... & \, 0
  \end{vmatrix}, 
 $$
 $x_{ij}$ are defined in (1), and $\tau_1$ is the smallest root $t>1$ of $C_G(t)=0.$ 

Theorem  3.1 in \cite{musin19} yields the following theorem. 
 \begin{theorem} \label{Th2} Let $G$ be a simple graph with $n$ vertices and  ${\mathcal R}(G)<\infty$.  Then $G$ can be embedded in $\mathbb R^d$  as a spherical two--distance set $S$ with $d=n-\mu(G)-1$ and $B(S)=\theta(G)$. If  
   $$\theta(G)+\mu(G)> n$$ then $S$ is a spherical counterexample for Borsuk's conjecture.  
\end{theorem}

%---------------------------------------------------------------------

\section{Remarks} 

%---------------------------------------------------------------------

\subsection{Bondarenko's counterexamples for Borsuk's conjecture.}
Bondarenko \cite{Bond14} found two spherical two--distance sets:\\
 $S_1=\{p_1,...,p_{416}\}$ in $\mathbb R^{65}$ with $\langle p_i, p_j \rangle =1/5$ or $\langle p_i, p_j \rangle =-1/15$ for $i\ne j$, $B(S_1)\ge 84$, \\ and  \\
 $S_2=\{p_1,...,p_{31671}\}$ in $\mathbb R^{782}$ with $\langle p_i, p_j \rangle =1/10$ or $\langle p_i, p_j \rangle =-1/80$ for $i\ne j$, $B(S_2)\ge 1377$. \\ 
 Jenrich's counterexample  in \cite{Jen} is \\
  $S_3=\{p_1,...,p_{352}\}$ in $\mathbb R^{64}$ with  $B(S_3)\ge 71$. 

% \medskip 

Let $S$, $|S|=n$,  be a two--distance set in $\mathbb R^d$. Then, see \cite{bannai1983,blokhuis1984}, we have
$$
n\le\frac{(d+1)(d+2)}{2}. 
 $$ 
 Bondarenko and Jenrich's counterexamples are spherical and satisfy the assumption of Theorem 1 in \cite{musin2009}. That yields a better bound  
 $$
n\le c_2(d):=\frac{d(d+1)}{2}. % \eqno (3)
 $$
 However, 
 $$|S_1|=416 < c_2(65)=2145, \; |S_2|=31671<c_2(782)=306153, \; |S_3|=352<c_2(64)=2080.$$ 
 Thus, the cardinalities of $S_k$ are much less than the maximum possible in their dimensions.

%%%%%%%%%%%%%%%%%%%%%%%%%%

\subsection{A possible strategy for finding counterexamples.}

Let $G$ be a graph with $n$ vertices. Then, see \cite{ES, musin19},% considered Euclidean representations of graphs. They  proved that
 
 \medskip
 
\noindent (i)  {\em $\rp(G)=n-1$, i.e. $\mu(G)=0$,  if and only if $G$ is a disjoint union of cliques.}
  
%  \medskip

%\noindent Moreover,  %they  have shown that 

\medskip

\noindent (ii) {\em If $d=\rp(G)\le n-2$, then a Euclidean representation of $G$ in ${\Bbb R}^d$ as a two--distance set with fixed minimum distance $a=1$  is uniquely defined up to isometry.}

\medskip

Let $\cl(G)$ denote the minimal disjoint union of cliques in $G$. If $C_0$ is a disjoint union of $m$ cliques $\{c_1,...,c_m\}$  and $\cl(G)=C_0$, then $\theta(G)=m$. Therefore, if 
$$m+\mu(G)>n $$
 then  $G$  is a counterexample for Borsuk's conjecture.

 Let us fix $C_0$ and consider $G$ for which $\cl(G)=C_0$.  Our goal is to find $G$ with largest possible $\mu$.  If $G=C_0$ then  by (i)  $\mu(G)=0$. In this case, adding any edge between clicks increases $\mu$. We can continue adding and removing edges between cliques, increasing $\mu$ as long as possible. If at some point we reach a local maximum of $\mu$, we can try to consider other graphs and find graphs with larger $\mu$  for them by adding and removing edges. Here, $C_0$ play an important role. In Bondarenko's counterexamples, almost all cliques have the same size. Perhaps such $C_0$ should be used to find counterexamples.

%Let us denote by $\Gamma(C_0)$ a graph whose  $\mu$  is maximal and the disjoint union of cliques is $C_0$. Then the clique covering number $\theta$ of this graph is $m$ and $\mu$ is the  for a given $C_0$. Let $$D_2(C_0):=m+\mu(\Gamma(C_0)).$$ If  $D_2(C_0)>n$ then  $\Gamma(C_0)$  is a counterexample for Borsuk's conjecture.  

%%%%%%%%%%%%%%%%%%%%%%%%%%%%%%%%%%%%%%%%%

\subsection{On the generalization of Theorem 1 for  $s$--distance sets}
First consider an equivalent definition of graph representations.  There are two Euclidean representation numbers $\rp(G)$ and $\rp(\bar G)$, where $\bar G$ is the graph complement of $G$. These numbers can be different. Indeed, there is an obvious relation between polynomials $C_G(t)$ and $C_{\bar G}(t)$. Namely,  $C_{\bar G}(t)$ is the reciprocal polynomial of $C_G(t)$. If $G$ or $\bar G$ is not the complete multipartite graph,  then   $\tau_0(G):=1/\tau_1(\bar G)$ is a zero of $C_G(t)$ and there are no more zeros in the interval $I:=[\tau_0(G),\tau_1(G)]$. Moreover, a two--distance set $S$ with distances 1 and $\sqrt{t}$ is well--defined only if $t\in I$ \cite{ES}. 

Let $G=(V(G),E(G))$ be a graph on $n$ vertices. We have $E(K_n)=E(G)\cup E(\bar G)$. Then it is can be considered as a coloring of $E(K_n)$ in two colors. 

%Hence $E(K_n)=E(G)\cup E(\bar G).$ %, \; \mbox{ where } \; E(G_1)\cap E(G_2)=\emptyset.$$ 
%Clearly, $G$ is uniquely defined by the equation $E(G)=E_1$. 

Let $L(e):=i$ if $e\in E_i$. Then  $L:E(K_n)\to \{1, 2\}$ is a coloring of $E(K_n)$.  A representation $L$  as a two--distance set is  an embedding $f$ of $V(K_n)$ into  ${\Bbb R}^d$ such that $\dist(f(u),f(v)))=a_i$ for $[uv]\in E_i$. For $s$--distance sets, it is convenient to swap  $G$ and $\bar G$  when the maximum distance is equal to 1, i.e. $a_1=1> a_2>0$.

\medskip

This definition can be extended to any number of colors. Let $L:E(K_n)\to \{1,\ldots,s\}$ be a coloring of the set of edges of a complete graph $K_n$. Then 
$$
E(K_n)=E(G_1)\cup\ldots\cup E(G_s), \; E(G_i):=\{e\in E(K_n): L(e)=i\}. 
$$
We say that  an embedding $f$ of the vertex set of $K_n$ into  ${\Bbb R}^d$  is a {\em Euclidean representation of a coloring $L$  in ${\Bbb R}^d$ as an $s$--distance set} if there are $s$ positive real numbers $a_1=1> \ldots> a_s$ such that $\dist(f(u),f(v)))=a_i$ if and only if $[uv]\in E(G_i)$.

Let  $\rp(L)$ denotes the smallest dimension of the embedding of $L$ as an $s$--distance set. Then a coloring $L$ defines a counterexample for Borsuk's conjecture only if 
$$
\theta(\bar G_1)>\rp(L)+1.
$$

It is easy to extend the definitions of  $A_G(t)$ and $C_G(t)$ in (1) for $s$--distance sets. 
In this case we have a matrix $A_L(t_2,\ldots,t_{s})$ and a polynomial $C_L(t_2,\ldots,t_{s})$, where $a_1=1$ and  $t_i=a_i^2<1$\, for $i=2,\ldots,s$. Similarly to $s=2$, for the Euclidean representation $L$ can be defined $M_L(t_2,\ldots,t_{s})$. Then this embedding  is spherical only if $H_L(t_2,\ldots,t_{s})$ is well defined, where 
$$H_L(t_2,\ldots,t_{s}):=\frac{M_L(t_2,\ldots,t_{s})}{C_L(t_2,\ldots,t_{s})}.$$

The main problem here is to develop an efficient algorithm that finds $t_2,...,t_s$ such that the rank of $A_L$, i.e. $\rp(L)$,  is  minimal. It is possible that a similar approach as in 3.2 can be used to find counterexamples  among $s$--distance sets to Borsuk's conjecture.

%We think that the Einhorn--Schoenberg theorem and several results from this paper can be generalized for representations of colorings $L$ as $s$--distance sets. 

%%%%%%%%%%%%%%%%%%%%%%%%%%%

\medskip

%\noindent{\bf Acknowledgment.}  I am very grateful to the reviewers of this paper for their great help in improving the text and helpful comments.

\medskip

\medskip

\medskip

\medskip

 O. R. Musin, School of Mathematical and Statistical Sciences, University of Texas Rio Grande Valley,  One West University Boulevard, Brownsville, TX, 78520.

 {\it E-mail address:} oleg.musin@utrgv.edu

\end{document}